\documentclass[3p,preprint]{elsarticle}

\usepackage{amsmath}
\usepackage{amsfonts}
\usepackage{amssymb,latexsym}
\usepackage{mathrsfs}
\usepackage{amsthm}
\usepackage{graphicx}

\newtheorem{theorem}{Theorem}[section]

\newtheorem{Remark}{Remark}[section]

\usepackage{color}
\usepackage{newlfont}
\usepackage{subfigure}
\usepackage{verbatim}
\usepackage{rotating}
\usepackage{multirow}
\usepackage{units}
\usepackage{bm}
\usepackage[english]{babel}
\usepackage[utf8]{inputenc}
\usepackage{algorithm}
\usepackage[noend]{algpseudocode}
\usepackage{booktabs}
\usepackage{caption}
\usepackage{hyperref}

\journal{Applied Numerical Mathematics}

\begin{document}

\begin{frontmatter}

\title{TetraFreeQ: tetrahedra-free quadrature on polyhedral elements}



\author[address-PD,address-GNCS]{Alvise Sommariva}
\ead{alvise@math.unipd.it}

\author[address-PD,address-GNCS]{Marco Vianello}
\ead{marcov@math.unipd.it}

\cortext[corrauthor]{Corresponding author}

\address[address-PD]{University of Padova, Italy}

\address[address-GNCS]{Member of the INdAM Research group GNCS}

\begin{abstract}
In this paper we provide a tetrahedra-free algorithm to compute low-cardinality  
quadrature rules 
with a given degree of polynomial exactness, positive weights and interior nodes 
on a polyhedral element with arbitrary shape. The key tools are the notion of 
Tchakaloff discretization set and the solution of moment-matching equations by Lawson-Hanson 
iterations 
for NonNegative Least-Squares. Several numerical tests are presented.  The method is implemented in Matlab as open-source software.
\end{abstract}

\begin{keyword} 
Polyhedral elements, tetrahedra-free algebraic quadrature, 
Positive Interior rules,  Tchakaloff theorem, Tchakaloff sets, Davis-Wilhelmsen theorem, NonNegative Least Squares.
\MSC[2020]65D32.
\end{keyword}


\end{frontmatter}

\section{Introduction}

Let $\Omega \subset {\mathbb{R}}^3$ is a polyhedron 
(either convex or nonconvex or even multiply-connected), and suppose that
for a given continuous function $f \in C(\Omega)$ one needs to approximate
its integral over $\Omega$ by a quadrature rule with nodes
$\{Q_j\}_{j=1,\ldots,\nu}\subset \Omega$ where $Q_j=(x_j,y_j,z_j)$ and
positive weights $\{w_j\}_{j=1,\ldots,\nu}$,
\begin{equation}\label{qformula}
I_{\Omega}(f)=\int_{\Omega} f(x,y,z) \,\,\, dx \, dy \, dz
\approx I_n(f)=\sum_{j=1}^\nu w_jf(Q_j)\;,
\end{equation}
that has Algebraic Degree of Exactness $ADE=n$, i.e. $I_{\Omega}(p)
= I_n(p)$ for all $p \in {\mathbb{P}}^3_n$ (the space of trivariate
polynomials of total degree not exceeding $n$). Such a kind of rules 
are often called PI (Positive Interior) in the literature 
and are optimally
stable by the positivity of the weights, being the quadrature
conditioning $\sum_j |w_j|/\left|\sum_{j}w_j\right|=1$.
As a relevant motivation, we may recall that the algebraic quadrature problem on polyhedra has been extensively 
studied in the literature, 
especially during the last decade, being at the core of many discretization 
methods for PDEs that use polyhedral meshes, such as polyhedral FEM and  
discontinuous Galerkin, as well as the more recent VEM. We may quote among 
others, 
without any pretence of exhaustivity, the papers 
\cite{AHP18,BDR17,MS11,SDAW14,SSVW17} with the numerous references therein.     

Despite of such manifest interest in the framework of numerical PDEs, open-source numerical codes to compute PI algebraic quadrature 
rules on polyhedra, 
in particular Matlab codes, do not seem to be readily available. 
On the other hand, differently from other approaches that do not provide an algebraic quadrature formula but a numerical integration algorithm, the availability of low-cardinality quadrature formulas allows to compute efficiently integrals involving forcing terms, without the need of an explicit polynomial approximation of such terms (cf. \cite[Rem. 6]{AHP18}). Moreover, since the polyhedral quadrature formulas can be computed once and for all for a fixed polyhedral mesh, 
even more flexibility is gained when different forcing terms 
have to be tested, or more generally when parametrized PDEs have to be solved, for example by reduced basis methods or other model order reduction techniques (cf., e.g., \cite{HRS16}).

Indeed, 
the main purpose of the present work is to provide an algorithm implemented by an open-source Matlab package, named {\sc TetraFreeQ}, for the computation of tetrahedra-free PI algebraic quadrature rules on arbitrary polyhedra 
(either convex or nonconvex or even multiply-connected), 
with cardinality $\nu\leq N$ where 
\begin{equation} \label{dim}
N=N_n=dim(\mathbb{P}_n^3)=(n+1)(n+2)(n+3)/6\;.
\end{equation} 
Existence and computability of such rules are guaranteed by the well-known Tchakaloff theorem \cite{T57} and the less known Wilhelmsen theorem on ``Tchakaloff sets'' \cite{W76}, the latter being at the core of our approach as described in the nest section.

\section{Tetrahedra-free quadrature}
In this section we show how to determine a quadrature rule with $ADE=n$, positive weights and interior nodes, on a general polyhedron $\Omega$, without needing a decomposition of the domain by means of nonoverlapping tetrahedra (tetrahedra-free quadrature, for short). It should be recalled that avoiding ``sub-tessellations'' for integration on polyhedral elements is a common 
requirement in the numerical PDEs framework, especially with high-order methods. 

In the sequel the polyhedron is assumed to be given via its polygonal faces, say $\{\mathcal{F}_i\}$, each with counterclockwise oriented vertices w.r.t. the outward normal.
The procedure works essentially as follows:
\vskip0.5cm
\noindent{\bf Algorithm TetraFreeQ} (tetrahedra-free quadrature on a general polyhedron $\Omega$)
\begin{itemize}
\item[$(i)$] compute the moment array 
\begin{equation} \label{moms}
\gamma=\{\gamma_k\}=\left\{\int_\Omega{\phi_k(x,y,z)\,dx\,dy\,dz}\right\}
\end{equation}
of a certain polynomial basis 
$\{\phi_k\}$ of ${\mathbb{P}}_n^3$, $k=1,\ldots,N$, by the divergence theorem and bivariate  
quadrature rules with ADE equal to $n+1$ on the polyhedron faces $\{{\cal{F}}_i\}_{i=1,\ldots,M}$
\item[$(ii)$] using an {\it{in-polyhedron}} routine \cite{H15}, in view of Wilhelmsen theorem \cite{W76}, determine iteratively a sufficiently dense pointset $\{P_l\}_{l=1,\ldots,L}$ inside $\Omega$ (a so-called ``Tchakaloff set'') such that the underdetermined moment-matching system 
\begin{equation} \label{mm}
V^t u=\gamma\;,\;\; \mbox{where}\;\; V=[v_{lk}]=
[\phi_k(P_l)]\in \mathbb{R}^{L\times N}
\end{equation}
has a nonnegative solution $u$ with at most $N\leq L$ positive components, computed via Lawson-Hanson NNLS algorithm \cite{LH95} applied to 
\begin{equation} \label{nnls}
\min_{u\geq 0}{\|V^tu-\gamma\|_2}
\end{equation}
\item[$(iii)$] once $\{P_l\}_{l=1,\ldots,L}$ and the weight vector $u$ are determined, the nonzeros of $u$ finally select the nodes of a PI-type quadrature rule with cardinality at most $N$.
\end{itemize}
\vskip0.5cm

In spite of the simplicity of this approach, there are many theoretical and practical aspects that deserve explanations. We stress that the key tools are Wilhelmsen theorem on the existence of ``Tchakaloff sets'' within 
sufficiently dense discrete subsets (a relevant but somehow overlooked result of multivariate quadrature theory), and the capability of Lawson-Hanson active-set iterative algorithm to compute a sparse solution to the NonNegative Least-Squares problem (made more efficient by the recent Deviation Maximization approach to QR factorization \cite{DM22,DMV20}). 

\subsection{Moment computation}
For an implementation of the algorithm sketched above, we intend to compute first the moments on the polyhedron $\Omega$ of the product Chebyshev basis of total degree $n$, relative to the smallest parallelepiped ${\cal{R}}=[a_1,b_1] \times [a_2,b_2] \times [a_3,b_3]$ containing $\Omega$.

In other words, defined the scaled Chebyshev polynomial in $[a,b]$
\begin{equation}\label{mc0}
{\tilde{T}}^{(a,b)}_{m}(t) := T_m\left(  \frac{2}{b-a} \cdot   \left(t-\frac{a+b}{2} \right)  \right)
\end{equation}
where $T_m$ is the Chebyshev polynomial of first kind, of degree $m$,  we compute the value of the moments
$$
\gamma_{\alpha}=\int_{\Omega} {\tilde{T}}^{(a_1,b_1)}_{\alpha_1}(x)  {\tilde{T}}^{(a_2,b_2)}_{\alpha_2}(y)  {\tilde{T}}^{(a_3,b_3)}_{\alpha_3}(z) {\mbox{ }} dx \, dy \, dz
$$
for all the triples $\alpha=(\alpha_1,\alpha_2,\alpha_3)$ with $\alpha_1,\alpha_2,\alpha_3 \in \mathbb{N}$ and $0\leq \alpha_1+\alpha_2+\alpha_3 \leq n$.

To this purpose, we recall that in view of the divergence theorem applied to a continuously differentiable vector field ${\mathbf{F}}=(f,g,h)$ defined on a neighborhood of $\Omega$, we have
\begin{equation}\label{mc1}
\int_{\Omega} \nabla \cdot {\mathbf{F}} {\mbox{ }} dx \, dy \, dz = \int_{\partial \Omega} {\mathbf{F}} \cdot {\mathbf{n}} {\mbox{ }} dS
\end{equation}
where ${\mathbf{n}}=(n_1,n_2,n_3)$ is the outward pointing unit normal at each point on the boundary ${\partial \Omega}$ and $ \nabla \cdot {\mathbf{F}} = {\frac{\partial f}{\partial x}}+{\frac{\partial g}{\partial y}}+{\frac{\partial h}{\partial z}}$ is the divergence of $\mathbf{F}$.
As a particular case of (\ref{mc1}), we consider ${\mathbf{F}}=(f,0,0)$ obtaining 
\begin{equation}\label{mc2}
\int_{\Omega} {\frac{\partial f(x,y,z)}{\partial x}} dx \, dy \, dz = \int_{\partial \Omega} n_1(x,y,z)  \, {{f(x,y,z)}}  {\mbox{ }} dS\;,
\end{equation}
and we apply (\ref{mc2}) to a certain $f$ such that $${\frac{\partial f(x,y,z)}{\partial x}}={\tilde{T}}^{(a_1,b_1)}_{\alpha_1}(x)  {\tilde{T}}^{(a_2,b_2)}_{\alpha_2}(y)  {\tilde{T}}^{(a_3,b_3)}_{\alpha_3}(z).$$
In this regard, we exploit a fundamental property of Chebyshev polynomials of first kind, namely that for $m \geq 2$
\begin{equation}\label{mc3}
\int_0^x T_m(t) dt = \frac{T_{m+1}(x)}{2(m+1)}- \frac{T_{m-1}(x)}{2(m-1)}.
\end{equation}
Consequently, setting $s= \frac{2}{b-a} \cdot (t-\frac{a+b}{2})$, by (\ref{mc0}) and (\ref{mc3}), for $m \geq 2$,
\begin{eqnarray}
\int_0^x {\tilde{T}}^{(a,b)}_{m}(t) dt &=& \int_0^x T_m \left(  \frac{2}{b-a}  \cdot  \left(t-\frac{a+b}{2} \right)  \right) dt \nonumber \\
&=& \int_{-\frac{a+b}{b-a}}^{\frac{2}{b-a} (x-\frac{a+b}{2}) } T_m(s) \frac{2}{b-a} ds \nonumber \\
&=& \frac{2}{b-a} \left( \int_{-\frac{a+b}{b-a}}^{0} T_m(s) ds  + \int_{0}^{\frac{2}{b-a} (x-\frac{a+b}{2})  } T_m(s) ds \right) \nonumber \\
&=& \frac{2}{b-a} \left( \int_{-\frac{a+b}{b-a}}^{0} T_m(s) ds  +  \frac{T_{m+1}( \frac{2}{b-a} (x-\frac{a+b}{2}) )}{2(m+1)}  -  \frac{T_{m-1}(\frac{2}{b-a} (x-\frac{a+b}{2}) )}{2(m-1)} \right) \nonumber \\
&=& \frac{2}{b-a} \left( \int_{-\frac{a+b}{b-a}}^{0} T_m(s) ds  +  \frac{ {\tilde{T}}^{(a,b)}_{m+1}(x)}{2(m+1)}  -  \frac{ {\tilde{T}}^{(a,b)}_{m-1}(x)}{2(m-1)} \right).
\end{eqnarray}
Thus, being the term $\frac{2}{b-a} \int_{-\frac{a+b}{b-a}}^{0} T_m(s) ds$ costant, setting 
$$
\mathcal{Z}^{(a,b)}_m(x)=\frac{1}{b-a} \left(  \frac{ {\tilde{T}}^{(a,b)}_{m+1}(x)}{m+1}  -  \frac{ {\tilde{T}}^{(a,b)}_{m-1}(x)}{m-1} \right)
$$
we get that for $m \geq 2$
\begin{eqnarray}\label{mc5}
{\tilde{T}}^{(a,b)}_{m}(x)=\frac{\partial } {\partial x}\mathcal{Z}^{(a,b)}_m(x).
\end{eqnarray}
It can be easily checked that (\ref{mc5}) holds also in the case $m=0$ and $m=1$, by defining respectively
$$
\mathcal{Z}^{(a,b)}_{0}(x)= x-\frac{a+b}{2}, \mbox{       } \mathcal{Z}^{(a,b)}_{1}(x)= \frac{1}{b-a} \left( x-\frac{a+b}{2} \right)^2,
$$
and that $\mathcal{Z}^{(a,b)}_m \in {\mathbb{P}}_{m+1}$.

Now, denoting by $\{{\cal{F}}_i\}_{i=1,\ldots,M}$ the polygonal faces that compose the boundary of the polyhedron $\Omega$, where  ${\mathbf{n}}^{(i)}=(n^{(i)}_1,n^{(i)}_2,n^{(i)}_3)$ is the outward normal of ${\cal{F}}_i$, $i=1,\ldots,M$, by the form of the divergence theorem (\ref{mc2}) for
$$
f_{\alpha}(x,y,z) = \mathcal{Z}^{(a_1,b_1)}_{\alpha_1}(x)  {\tilde{T}}^{(a_2,b_2)}_{\alpha_2}(y)  {\tilde{T}}^{(a_3,b_3)}_{\alpha_3}(z)
$$
from the application of (\ref{mc5}) with $a=a_1$ and $b=b_1$, we have
\begin{eqnarray}\label{mc6}
\gamma_{\alpha} &=& \int_{\Omega} {\tilde{T}}^{(a_1,b_1)}_{\alpha_1}(x)  {\tilde{T}}^{(a_2,b_2)}_{\alpha_2}(y)  {\tilde{T}}^{(a_3,b_3)}_{\alpha_3}(z) {\mbox{ }} dx \, dy \, dz \nonumber \\
&=& \int_{\Omega} \frac{\partial } {\partial x}\mathcal{Z}^{(a_1,b_1)}_{\alpha_1}(x) {\tilde{T}}^{(a_2,b_2)}_{\alpha_2}(y)  {\tilde{T}}^{(a_3,b_3)}_{\alpha_3}(z) {\mbox{ }} dx \, dy \, dz \nonumber \\
&=& \int_{\Omega} \frac{\partial } {\partial x} \left(\mathcal{Z}^{(a_1,b_1)}_{\alpha_1}(x) {\tilde{T}}^{(a_2,b_2)}_{\alpha_2}(y)  {\tilde{T}}^{(a_3,b_3)}_{\alpha_3}(z) \right) {\mbox{ }} dx \, dy \, dz \nonumber \\
&=& \int_{\partial \Omega} n_1(x,y,z) \mbox{  } \mathcal{Z}^{(a_1,b_1)}_{\alpha_1}(x) {\tilde{T}}^{(a_2,b_2)}_{\alpha_2}(y)  {\tilde{T}}^{(a_3,b_3)}_{\alpha_3}(z) {\mbox{ }} dS \nonumber \\
&=& \sum_{i=1}^M \int_{ {\cal{F}}_i} n^{(i)}_1 \mbox{  } \mathcal{Z}^{(a_1,b_1)}_{\alpha_1}(x) {\tilde{T}}^{(a_2,b_2)}_{\alpha_2}(y)  {\tilde{T}}^{(a_3,b_3)}_{\alpha_3}(z) {\mbox{ }} dS\;. \nonumber \\
\end{eqnarray}
Since the outward normal $n_1$ is constant on the face ${\cal{F}}_i$, so is $n^{(i)}_1$, and thus the integrand 
$$
 n^{(i)}_1 \mbox{  } \mathcal{Z}^{(a_1,b_1)}_{\alpha_1}(x) {\tilde{T}}^{(a_2,b_2)}_{\alpha_2}(y)  {\tilde{T}}^{(a_3,b_3)}_{\alpha_3}(z)
$$
is a polynomial of degree not exceeding $n+1$ on the $i$-th face, with $i=1,\ldots,M$, and consequently each moment $\gamma_{\alpha}$  can be recovered by 
a quadrature rule with ADE equal to $n+1$ for all the polygons ${\cal{F}}_1, \ldots, {\cal{F}}_M$.

Having this in mind, let ${\cal{F}} \subset \partial \Omega$ be a (polygonal) face of the polyhedron $\Omega$, and suppose that the vertices are oriented counterclockwise w.r.t. its outward normal ${\mathbf{n}}$. Since it is easier to determine a quadrature rule on a planar polygon ${\cal{F}}^*$, we apply a rototranslation $\Phi$, such that ${\cal{F}}^*=\Phi({\cal{F}})$ is a polygon belonging to the plane $\pi^*=\{(x,y,0), \,\, x,y \in \mathbb{R}\}$, oriented counterclockwise w.r.t. the vector $(0,0,1)$.

A quadrature formula on the planar polygon ${\cal{F}}^*$ with internal nodes $\{P_i\}$ and positive weights $\{w_i\}$, having ADE equal to $n+1$, can be determined as described in \cite{BSV19}, by applying a rule with ADE equal to $n+1$ on each element of a minimal triangulation of ${\cal{F}}^*$ or alternatively by a triangulation free routine as discussed in {\cite{SV21}}. Since the rototranslation $\Phi$ is an affine (bijective) map with Jacobian modulus equal to 1, we conclude that $\{\Phi^{-1}(P_i)\}$ and $\{w_i\}$ are the nodes and the weights of a quadrature rule over the polygonal face ${\cal{F}} \subset \partial \Omega$, with ADE equal to $n+1$ and consequently we are able to compute the set of moments $\{\gamma_{\alpha}\}$.

\begin{Remark} {\em The mapping $\phi$ can be easily obtained by applying a rotation $R$ to the plane $\pi$ containing the vertices of the polygon $\cal F$, so that $R \pi$ is parallel to the $xy$-plane $\pi^*=\{(x,y,0), \,\, x,y \in \mathbb{R}\}$ and then by applying a vertical translation $\tau$ so that $\hat{\cal F}:=R{\cal F}+\tau \subset \pi^*$. 
If a counterclockwise orientation of $\cal F$ is given by the user, then it is straightforward to obtain the same for $\hat{\cal F}$. For the quadrature routines over planar polygons it is a typical requirement that the orientation is counterclockwise. If the planar polygon $\hat{\cal F}$ is oriented clockwise w.r.t. the {\it{z}}-axis then, applying an additional rotation, we can take as ${\cal F}^*$ the polygon symmetric to $\hat{\cal F}$ w.r.t. the {\it{x}}-axis, otherwise it is sufficient to set ${\cal F}^*=\hat{\cal F}$. 
We recall that a necessary and sufficient condition for a polygon with vertices $(x_j,y_j,0)$, $j=1,\ldots,m$ to be oriented counterclockwise is that its signed area  
$\sum_{j=1}^{m-1} (x_{j+1}-x_j)(y_{j+1}+y_j)$ be negative 
(for further details see \cite{Weiss}).

In some instances, instead of the orientation of the vertices of ${\cal F}$, it is available its outward normal $\mathbf n$ (w.r.t. $\Omega$). In this case, the procedure is the similar to that of the previous remark. If in particular, after the application of the rotation $R$ (and the translation $\tau$), the normal to ${\cal {\hat{F}}}$ correspond to $\Phi({\mathbf{n}})=(0,0,-1)$ then as ${\cal F}^*$ we can consider the polygon symmetric to $\hat{\cal F}$ w.r.t. the x-axis, otherwise we set ${\cal F}^*=\hat{\cal F}$.
}
\end{Remark}

\subsection{Tchakaloff sets and computation of Tchakaloff-like rules}

In the previous subsection we have implemented step $(i)$ of the Algorithm TetraFreeQ sketched at the beginning 
of Section 2. With a little abuse of notation, we shall consider an ordering of the multi-indices $\alpha=(\alpha_1,\alpha_2,\alpha_3)$, $0\leq \alpha_1+\alpha_2+\alpha_3\leq n$, e.g. the graded lexicographical ordering, and we shall term 
$\{\phi_k(x,y,z)\}$ the product Chebyshev basis and $\{\gamma_k\}$ the corresponding moments, $1\leq k\leq N$.  

Now, step $(ii)$ of Algorithm TetraFreeQ rests on the following relevant (though often overlooked) result proved by Wilhelmsen in {\cite{W76}} where the key notion of {\em Tchakaloff set} is introduced (extending a result of Davis \cite{D67} to rather general functional spaces, including polynomial spaces).

{\vspace{0.1cm}}

\begin{theorem} (Wilhelmsen, 1976) Let $\Psi$ be the linear span of continuous, 
real-valued, linearly independent functions $\{\phi_k\}_{k=1,\ldots,N}$ 
defined on a compact set $\Omega \subset {\mathbb{R}}^d$. 
Assume that $\Psi$ satisfies the Krein condition (i.e. there is at least 
one $f \in \Psi$ which does not vanish on $\Omega$) and that $L$ is a positive 
linear functional, i.e. $Lf>0$ for every $f \geq 0$ not vanishing everywhere in $\Omega$. 

If $\{P_l\}_{i=1}^{\infty}$ is an everywhere dense subset of $\Omega$, 
then for sufficiently large $L$, the set $X=\{P_l\}_{l=1,\ldots,L}$ is 
a Tchakaloff set, i.e. there exist weights $w_j > 0$, $j=1,\ldots,\nu$, and nodes $\{Q_j\}_{j=1,\dots,\nu}
\subset X \subset \Omega$, 
with $\nu={\mbox{card}}(\{Q_j\}) \leq N$, such that 
\begin{equation} \label{tchcub}
Lf=\sum_{j=1}^\nu w_j f(Q_j)\;, {\mbox{   }} \forall f \in \Psi\;.
\end{equation}
\end{theorem}

\vskip1cm

In view of this existence result, we can the implement step $(ii)$ in our framework, that is $L f=I_\Omega(f)$ and $\Psi=\mathbb{P}_n^3$, by the following iterative substeps: 

\begin{itemize} 

\item[$(ii1)$] generate a set of uniform random or low-discrepancy points (e.g. Halton points) with cardinality $K$ in the smallest box 
$[a_1,b_1] \times [a_2,b_2] \times [a_3,b_3]$ containing $\Omega$

\item[$(ii2)$] determine those points belonging to $\Omega$ by an in-polyhedron routine like 
\cite{H15}, say $X=\{P_l\}_{l=1,\ldots,L}$ with $L\leq K$, and compute the Vandermonde-like matrix $V=[v_{lk}]=[\phi_k(P_l)]\in \mathbb{R}^{L\times N}$ 

\item[$(ii3)$] solve the NNLS problem $\min_{u\geq 0}{\|V^tu-\gamma\|_2}$ by Lawson-Hanson 
active-set method 

\item[$(ii4)$] if the residual $\|V^tu-\gamma\|_2<\varepsilon$ (where $\varepsilon$ is a given tolerance) \\$\mbox{}$ then set $\{w_j\}=\{u_l>0\}$ and $\{Q_j\}=\{P_l:\,u_l>0\}$  
\\$\mbox{}$ else goto $(ii1)$ increasing $K$, e.g. by a fixed factor $\theta>1$ 

\end{itemize}

We stress that 
since the set of random points becomes {\it{denser}} along the iterations, in view of Wilhelmsen theorem above termination of the iterative cycle is (at least 
theoretically) guaranteed. Substep $(ii3)$ is based on the fact that Lawson-Hanson algorithm 
computes a sparse solution of the NNLS problem, and can be performed in Matlab by resorting to 
the built-in routine {\tt lsqnonneg}. On the other hand, a substantial acceleration can be obtained via the recent Deviation-Maximization approach to $QR$ factorizations which are 
the core of Lawson-Hanson method; cf. \cite{DODM22,DM22,DMV20} for a detailed analysis of this approach. 

\begin{Remark} 
{\em 
({\bf Tetrahedra-based quadrature\/}) 

As already stressed, the aim of the present work is to provide a completely ``tetrahedra-free'' 
approach to quadrature on general polyhedral elements, implemented in Matlab. On the other hand, a classical approach 
for polyhedral quadrature, whenever a nonoverlapping ``tetrahedralization''  
of the polyhedron is available, 
consists in collecting and summing up the PI-formulas on the single tetrahedra by integration additivity.

A first problem consists in the computation of the tetrahedralization. 
If $\Omega$ is convex or star-shaped with a known center, this operation is 
straightforward, but the task becomes not trivial for a general polyhedron 
that may have a difficult geometry. We recall, incidentally, that 
differently from the 2D case of polygons, there exist polyhedra that 
cannot be triangulated using only their vertices (e.g. the well-known 
Sch\"{o}nardt polyhedron, cf. \cite{Scho}).   

Moreover, though Matlab provides a {\sl{minimal}} 
triangulation of a general polygonal region via the environment 
{\tt{polyshape}}, similar codes do not seem available for polyhedra. 
As an alternative, sometimes useful, the built-in {\tt{alphashape}} 
command creates a bounding volume that envelops a given 3D point cloud 
pointset. If the polyhedron of 
interest can be defined in such a way, then {\tt{alphaTriangulation}} 
returns a tetrahedralization that defines the domain of the alpha shape 
(see {\cite{MW}} for additional details). 

A second problem concerns the rule to be used in each tetrahedron. 
For degrees of exactness up to $n=20$, 
there 
are in literature several PI-type rules on the reference tetrahedron (the simplex) 
with vertices  $(1, 0, 0)$,  $(0, 1, 0)$, $(0, 0 ,0)$, $(0, 0, 1)$, which have near-minimal 
cardinality (see e.g. {\cite{JS20,SH}} and Table 1). For $n > 20$, one can 
use the 
well-established Stroud rule \cite[p.28-32]{STR}, again of PI-type,  
that in general is not minimal but still easy to be 
implemented and with a moderate cardinality equal to $\lceil {\frac{n+1}{2}} \rceil^3$. 
All these rules are typically written in barycentric coordinates $(\lambda_1,\lambda_2,\lambda_3,\lambda_4)$ and can be applied  
to a general tetrahedron with vertices $A,B,C,D$ by the standard change of coordinates $\lambda_1 A+\lambda_2 B+\lambda_3 C+\lambda_4 D$. 

\vskip0.5cm
\begin{table}[ht]
\begin{center}\footnotesize
{\renewcommand{\arraystretch}{1.1}
\begin{tabular}{||| c | c ||| c | c ||| c | c ||| c | c |||}
\hline
{\mbox{deg}} & {\mbox{card}} & {\mbox{deg}} & {\mbox{card}} & {\mbox{deg}} & {\mbox{card}} & {\mbox{deg}} & {\mbox{card}}    \\
\hline
$ 1 $& $ 1 $ & $6$ & $ 23$& $  11$ & $ 94 $ &$  16$ & $247 $\\
$2$ & $  4 $ & $7$ & $ 31 $ & $  12$ & $117 $ &$  17$ & $288 $\\
$3$ & $  6 $& $8$ & $ 44 $& $  13$ & $144 $ &$  18$ & $338$\\
$4$ & $ 11 $&$9$ & $ 57 $& $  14$ & $175 $&$  19$ & $390 $ \\
$5$ & $ 14 $&$  10$ & $ 74 $& $  15$ & $207 $&$  20$ & $448 $\\
\hline
\end{tabular}
        }
\caption{\small{Cardinality of near-minimal rules on the reference 
tetrahedron.}}
\label{tabT01}
\end{center}
\end{table}

Of course the overall cardinality of the resulting rule is proportional to the number 
of tetrahedra and can be much larger than $N=(n+1)(n+2)(n+3)/6=dim(\mathbb{P}_n^3)$. 
In these cases, 
one can resort 
to Tchakaloff-like compression of the corresponding discrete measure, 
obtaining a rule with at most $N$ re-weighted nodes that are a subset 
of the starting ones. Such a compression consists in computing 
a sparse positive solution of a moment-matching 
system like (\ref{mm}), where the points and 
moments are given by the high-cardinality rule itself. This approach inserts in the more general 
setting of ``Caratheodory-Tchakaloff discrete measure compression'', and can be implemented by both Linear and Quadratic Programming methods; cf. e.g. \cite{Hay21,PSV,SV15,T15} with the references therein. 

In our Matlab package we have added, for the purpose of comparison, the routines 
that provide a compressed Tchakaloff-like rule along the lines sketched above, starting from a tetrahedralization 
of the polyhedral element. We notice that even this approach does not 
seem to be yet readily available, at least in Matlab, within existing open-source packages, 
and thus could be of independent practical interest, especially because it resorts to the best tetrahedral rules available in the literature. 
\/}
\end{Remark}
 
\section{Numerical tests}

In this section we test Algorithm TetraFreeQ on several polyhedral domains. All the routines have been implemented in Matlab and are freely available at  {\cite{SS}}. The numerical experiments were made using Matlab R2022a, on an Apple MacBook Pro with M1 Chip and 16 GB of RAM.

\begin{figure}[!htbp]
  \centering
   {\includegraphics[scale=0.40,clip]{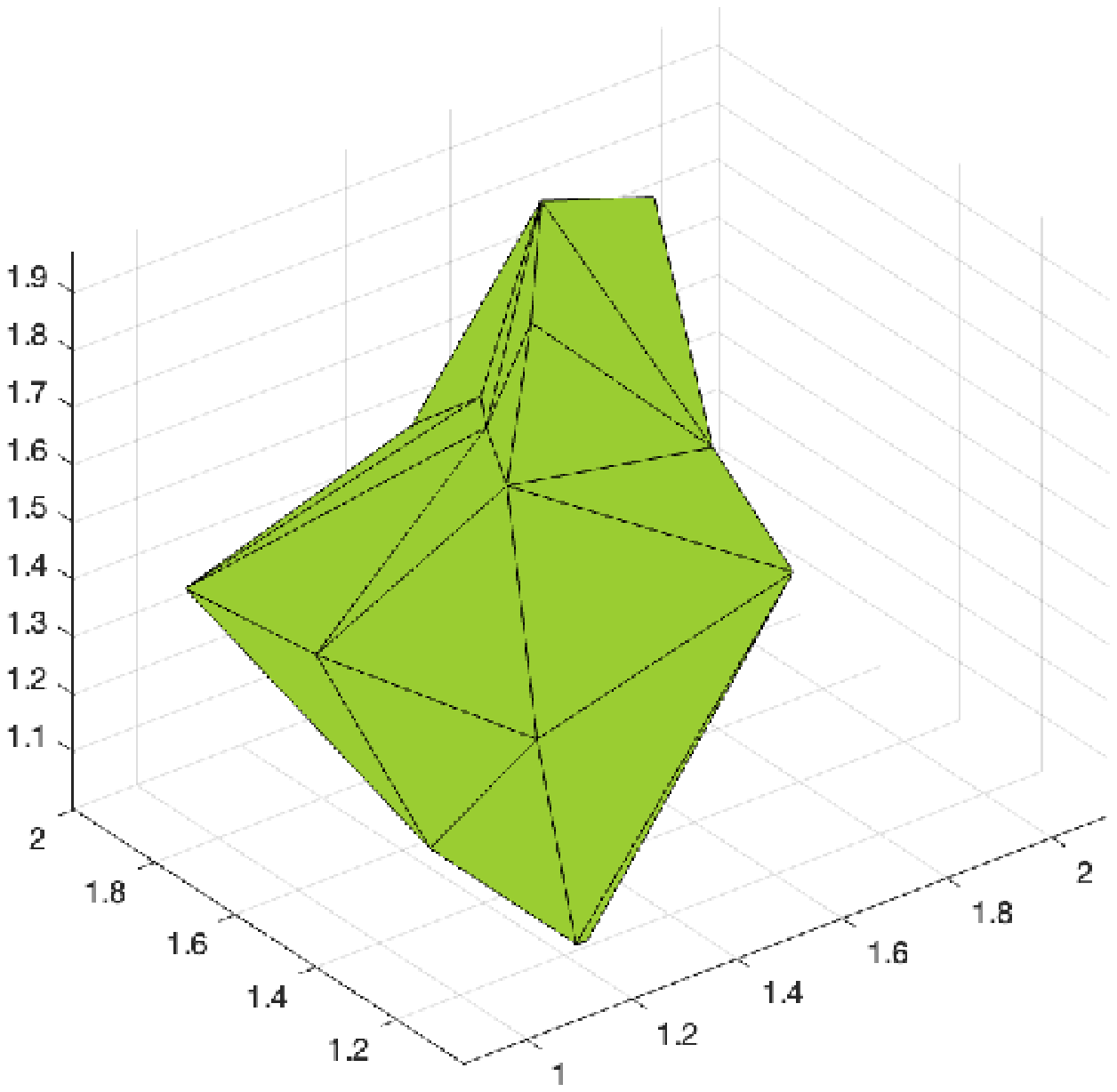}}
   {\includegraphics[scale=0.40,clip]{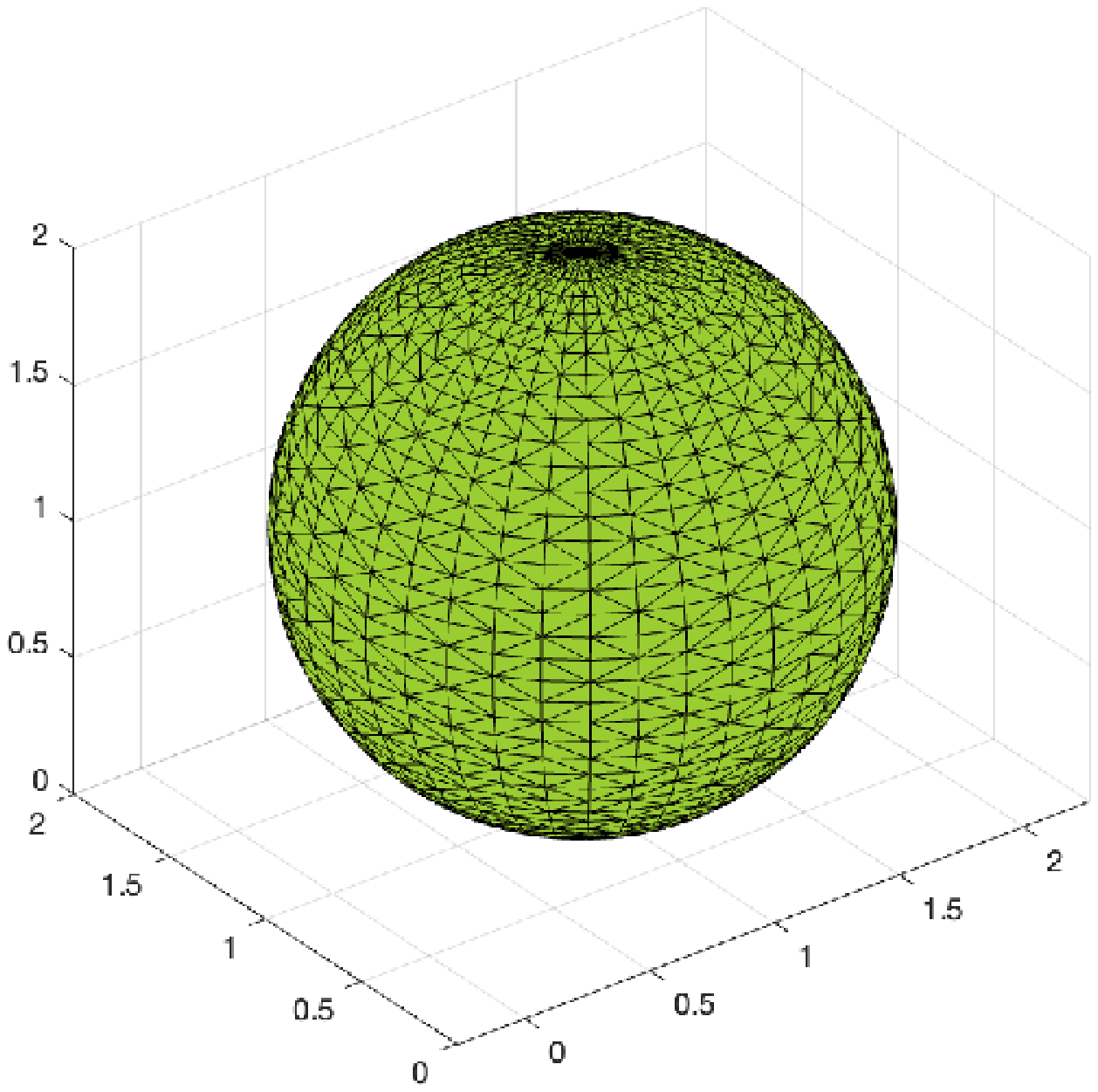}}
      {\includegraphics[scale=0.40 ,clip]{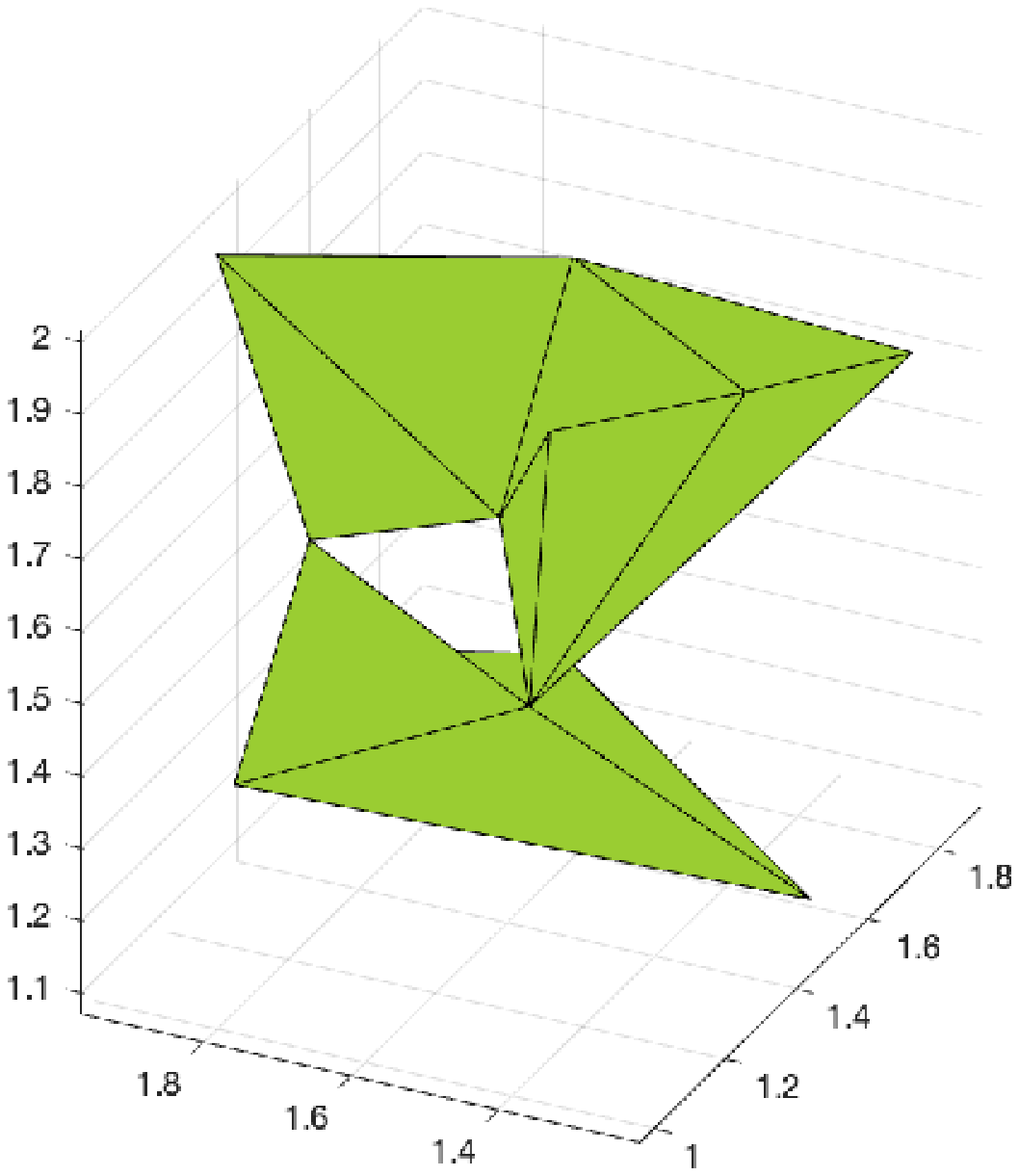}}
 \caption{Examples of polyhedral domains. {{Left:}} $\Omega_1$ (nonconvex, 20 facets); {{Center}}: $\Omega_2$ (convex, 760 facets); {{Right}}: $\Omega_3$ (multiply connected, 20 facets).}
 \label{fig_2CRI}
 \end{figure}

We consider three polyhedra $\Omega_j$, $j=1,2,3$ 
with triangular facets, where $\Omega_1$ is nonconvex with 30 facets, $\Omega_2$ is convex (an approximation of a sphere with 760 facets) and $\Omega_3$ is multiply-connected with 20 facets and a hole (see Figure \ref{fig_2CRI}). The domains, as well as the facets and the tetrahedralization (used for comparison), are obtained by Matlab built-in command  {\tt{alphashape}}, on suitable point clouds.

In particular we compute on the three polyhedra:
\begin{itemize}
\item a reference rule $Q_{T}$ corresponding 
to a tetrahedralization (cf. Remark 2.2) 

\item the corresponding compressed tetrahedral rule $Q_{TC}$

\item the tetrahedra-free rule $Q_{TF}$ produced by Algorithm TetraFreeQ. 

\end{itemize} 

In our numerical tests, in $(ii1)-(ii4)$ of Algorithm TetraFreeQ we have chosen 
the 
parameters $\varepsilon=10^{-14}$, $\theta=4$ and an initial value $K=\lceil 8Nvol(B)/vol(\Omega)\rceil$ where $B$ is the smallest bounding box for $\Omega$. 

In order to test numerically polynomial exactness of $Q_{TF}$, in Figure 3 we display for each domain the relative errors 
\begin{equation} \label{relerr}
E(g_k)=\frac{|Q_{T}(g_k)-Q_{TF}(g_k)|}{|Q_T(g_k)|}
\end{equation}
corresponding to  
100 trials of the polynomial 
\begin{equation} \label{random-p}
g_k(x,y,z)=(a_k+b_kx+c_ky)^n
\end{equation}
where $n$ is the algebraic degree of precision of the rule and $(a_k,b_k,c_k)\in [0,1]$ are uniform 
random coefficients. In the same figures, we have also represented with a circle the logarithmic averages $\sum_{k=1}^{100}{\log(E(g_k))/100}$. 

As additional investigation, we also compute the relative errors of these rules in integrating on $\Omega_j$, 
$j=1,2,3$, the three test functions 
with different degree of regularity 
\begin{equation} \label{testf}
f_1(x,y,z)=\exp(-\|P-P_0\|_2^2)\;, \;\;f_2(x,y,z)=\|P-P_0\|_2^5\;,\;\;
f_3(x,y,z)=\|P-P_0\|_2\;,
\end{equation}
where $P=(x,y,z)$ and $P_0(x_0,y_0,z_0)$ was suitably chosen in the domains, more precisely $P_0=(1.5,1.5,1.5)$ for $\Omega_1$ and $\Omega_3$, $P_0=(1,1,1)$ for $\Omega_2$; cf. Figures 4-6. Moreover, in Table 2 we report the cardinalities of the pointsets from which the tetrahedra-free rules are
extracted, varying the degree from 1 to 10.

In view of these numerical experiments, we can see that the rules $Q_{TC}$ and $Q_{TF}$ give errors of comparable size, whereas the errors 
of $Q_T$ are smaller. This is not surprising, since the basic tetrahedra-based rule uses a much larger number of nodes than $N$, which turns out to be the number of nodes of both $Q_{TC}$ and $Q_{TF}$, with a ratio increasing with the number of facets (being for example more than 100 with the 760-facets sphere-like polyhedron, cf. Figure 2). We stress once again, however, that $Q_{TF}$ avoids completely any tetrahedralization of the polyhedron, a desiderable feature in many applications (such as polyhedral FEM/VEM methods).

Concerning the cputime required by the moment computation in TetraFreeQ, it grows with the algebraic degree of precision and of course with the number of facets, ranging in our tests between $10^{-4}$ and $10^{-3}$ seconds. On the other hand, there is numerical evidence that the bottleneck of TetraFreeQ is due to the in-polyhedron routine and the NNLS algorithm, one dominating on the other depending on the algebraic degree of exactness and on the number of facets, with an overall cputime of the quadrature rule construction 
varying from $10^{-2}$ seconds for $n=1,2,3$ up to $10^{0}$ seconds for $n=10$.

\begin{figure}[!htbp]
  \centering
   {\includegraphics[scale=0.32,clip]{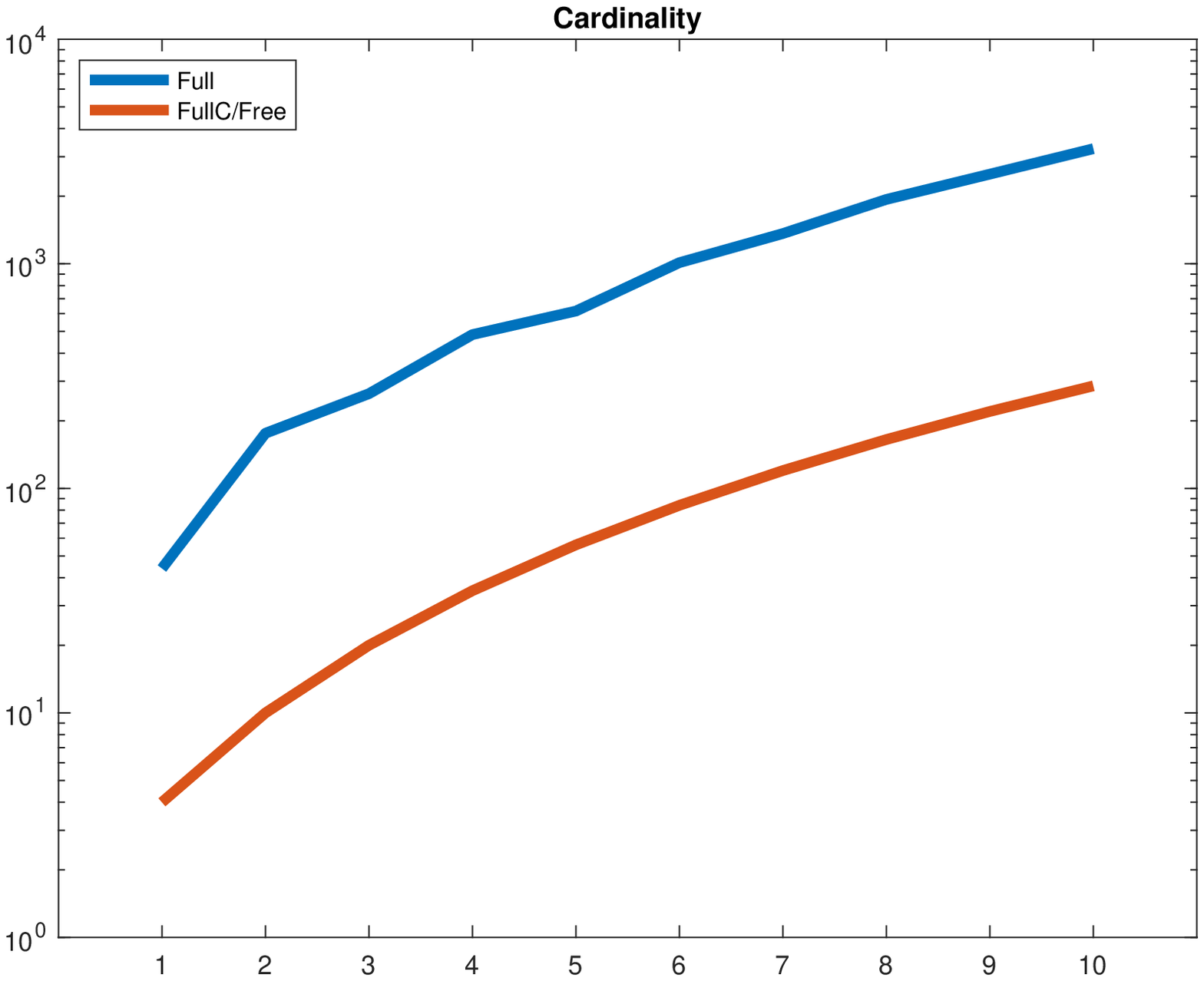}}
   {\includegraphics[scale=0.32,clip]{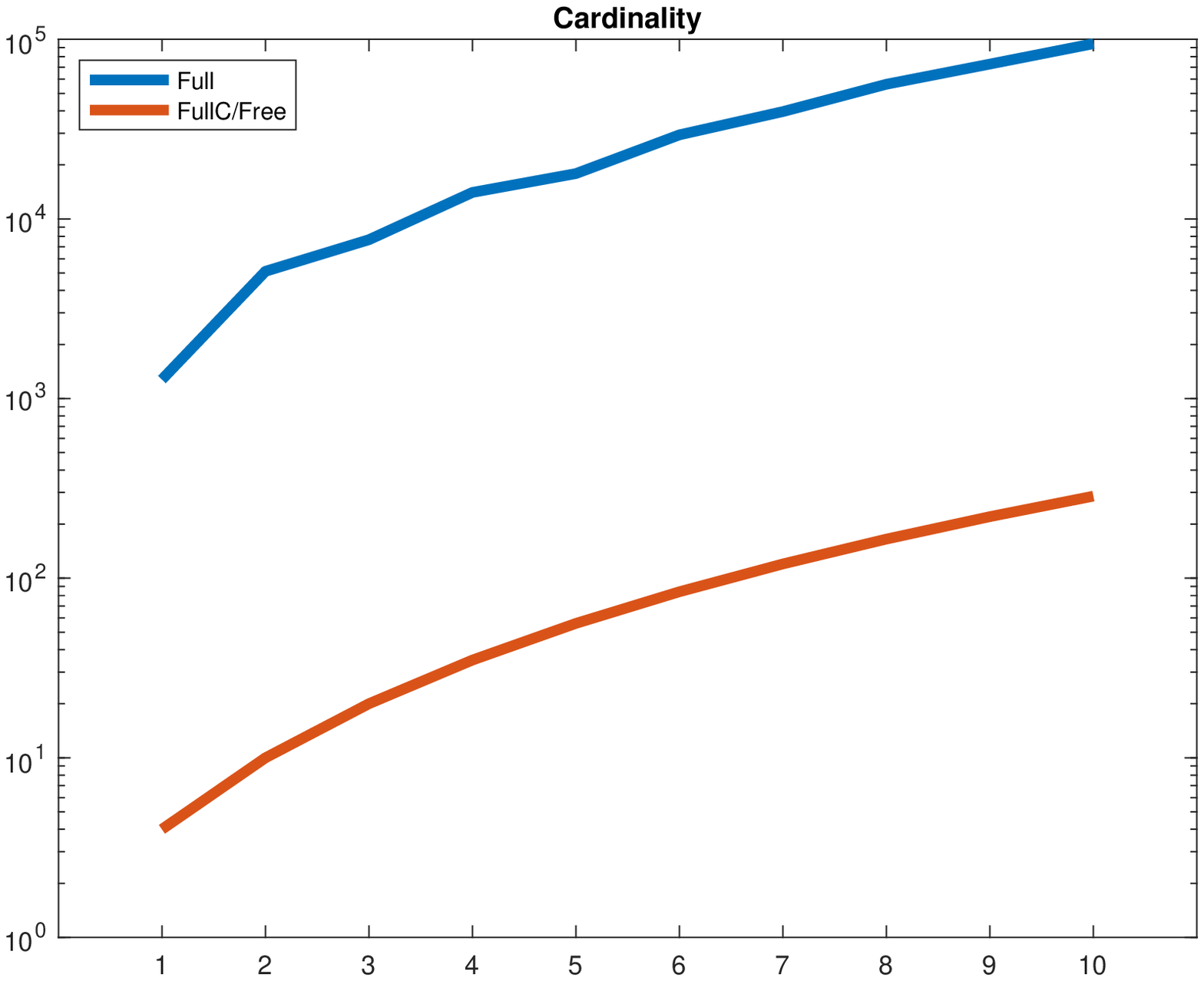}}
      {\includegraphics[scale=0.32,clip]{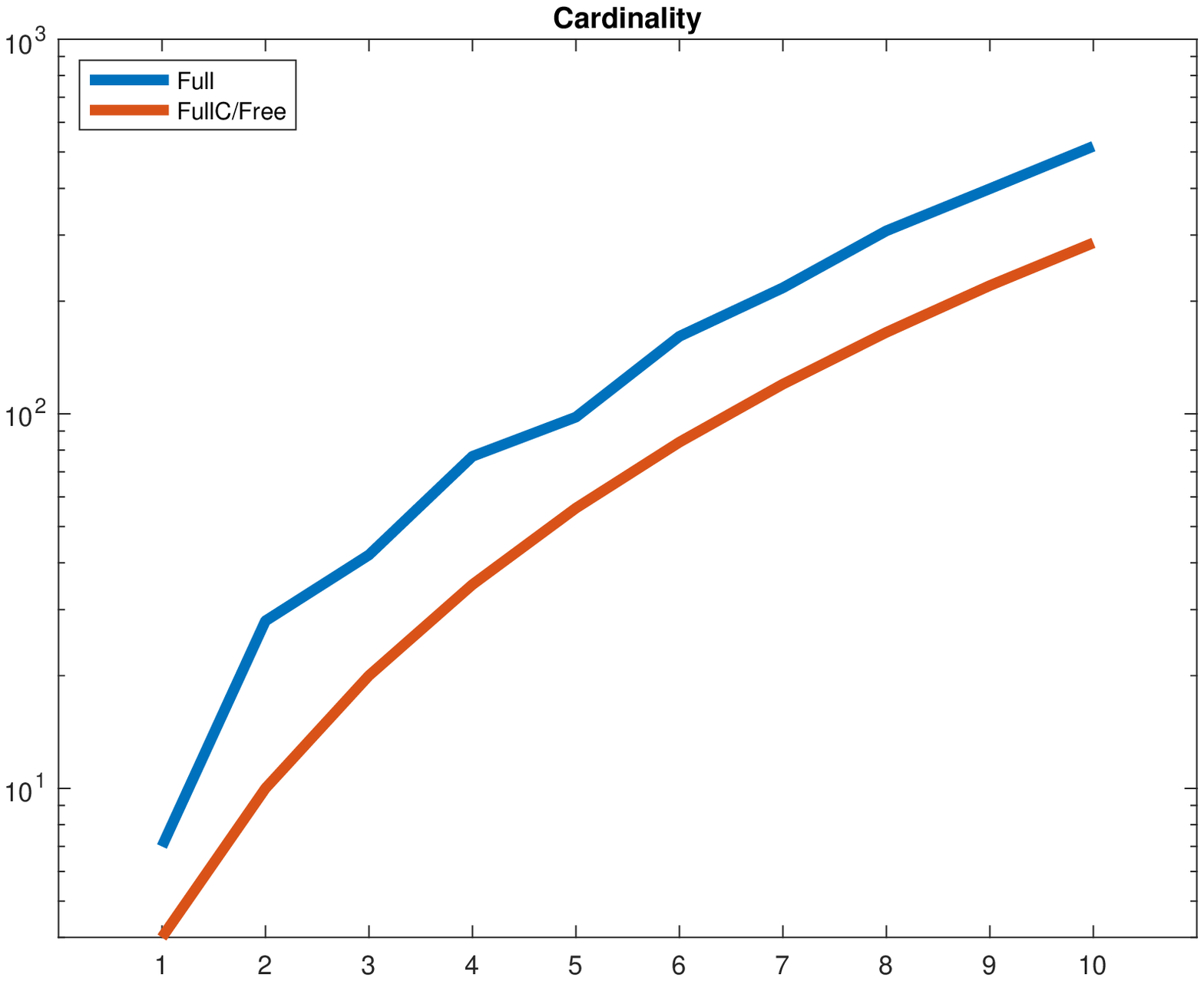}}
 \caption{Cardinality versus the exactness degree on the three polyhedra of Figure 1: uncompressed tetrahedra-based rule (blue), compressed tetrahedra-based and tetrahedra-free rule (brown).}
 \label{fig_3CRI}
 \end{figure}

\begin{figure}[!htbp]
  \centering
   {\includegraphics[scale=0.32,clip]{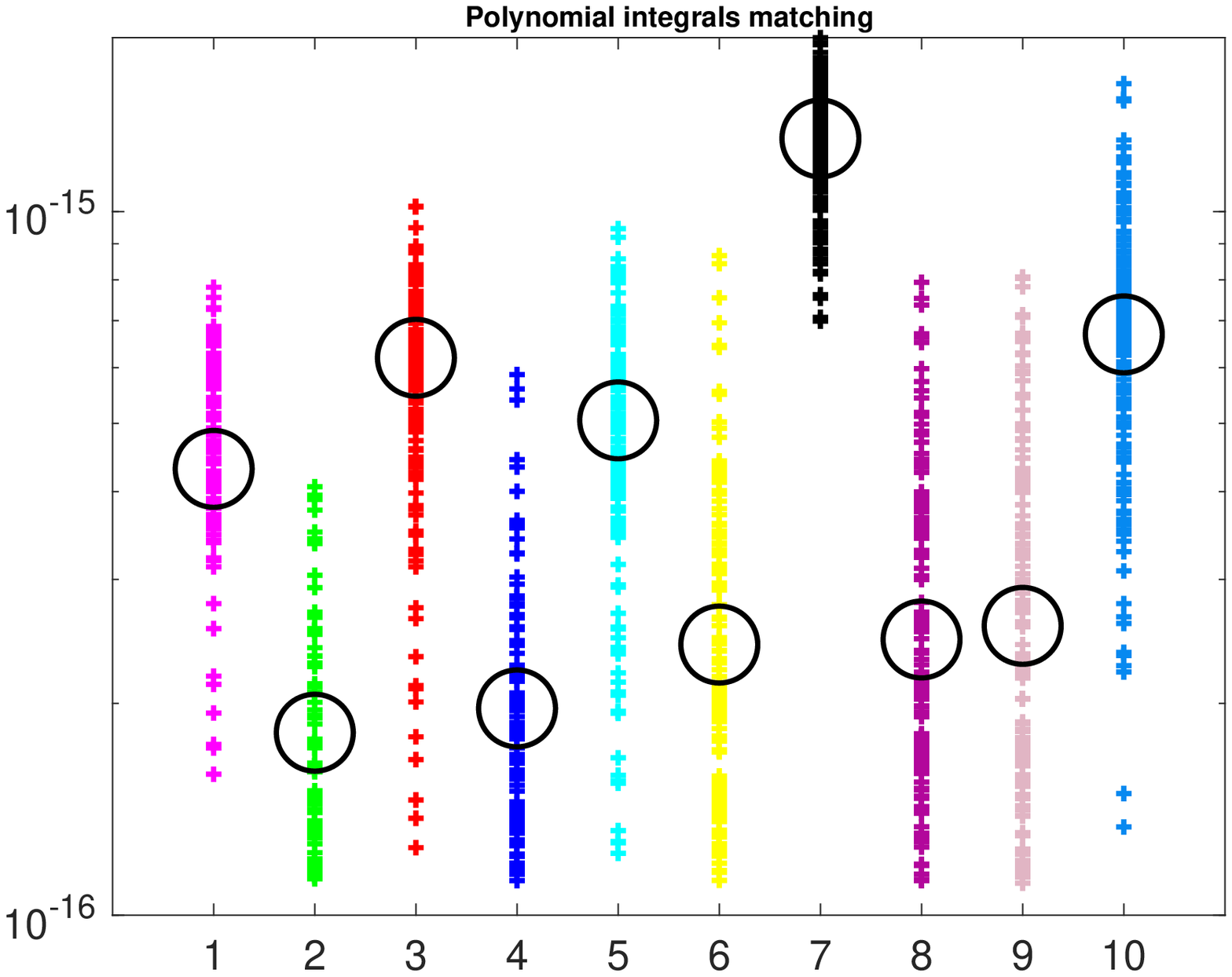}}
   {\includegraphics[scale=0.32,clip]{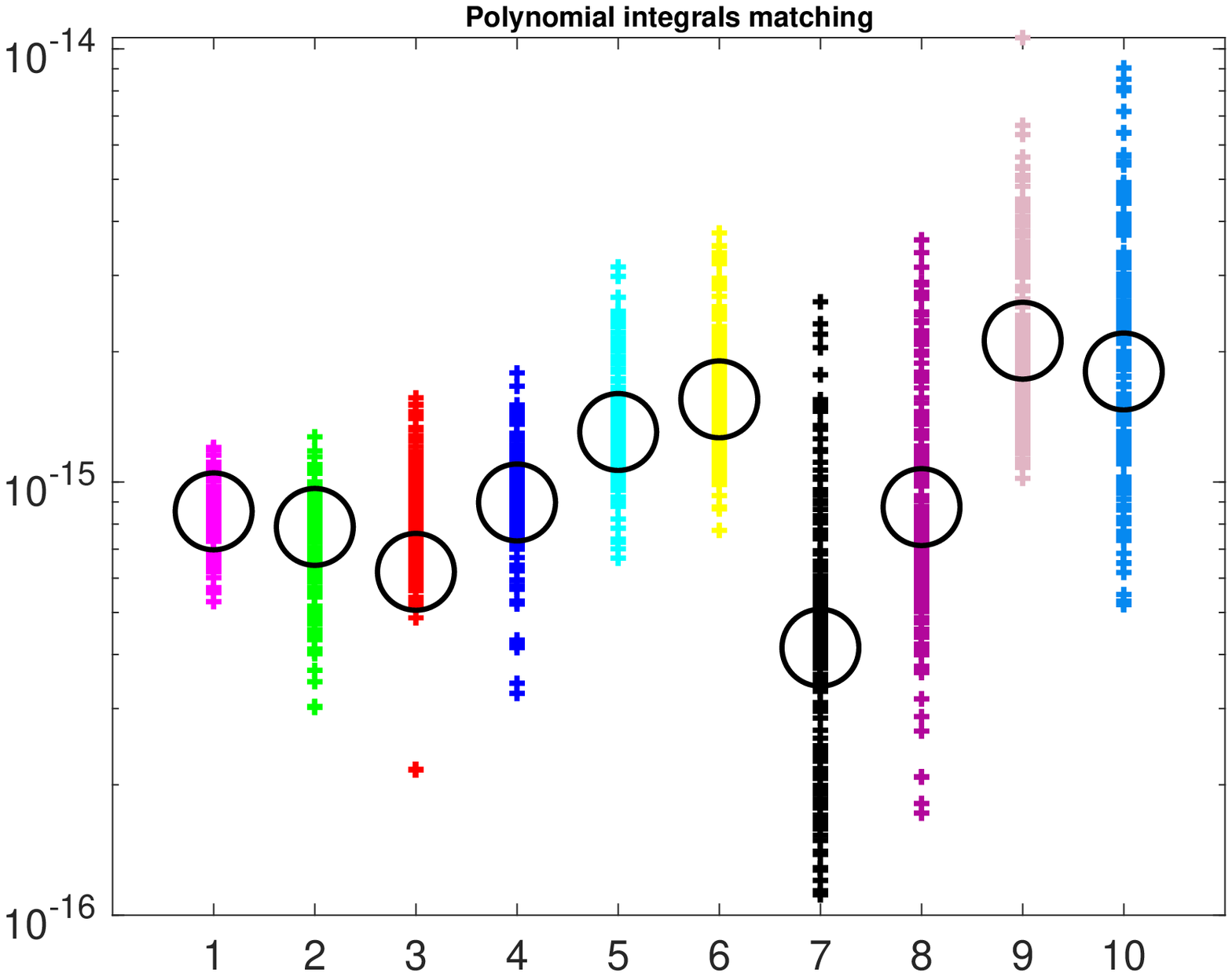}}
      {\includegraphics[scale=0.32,clip]{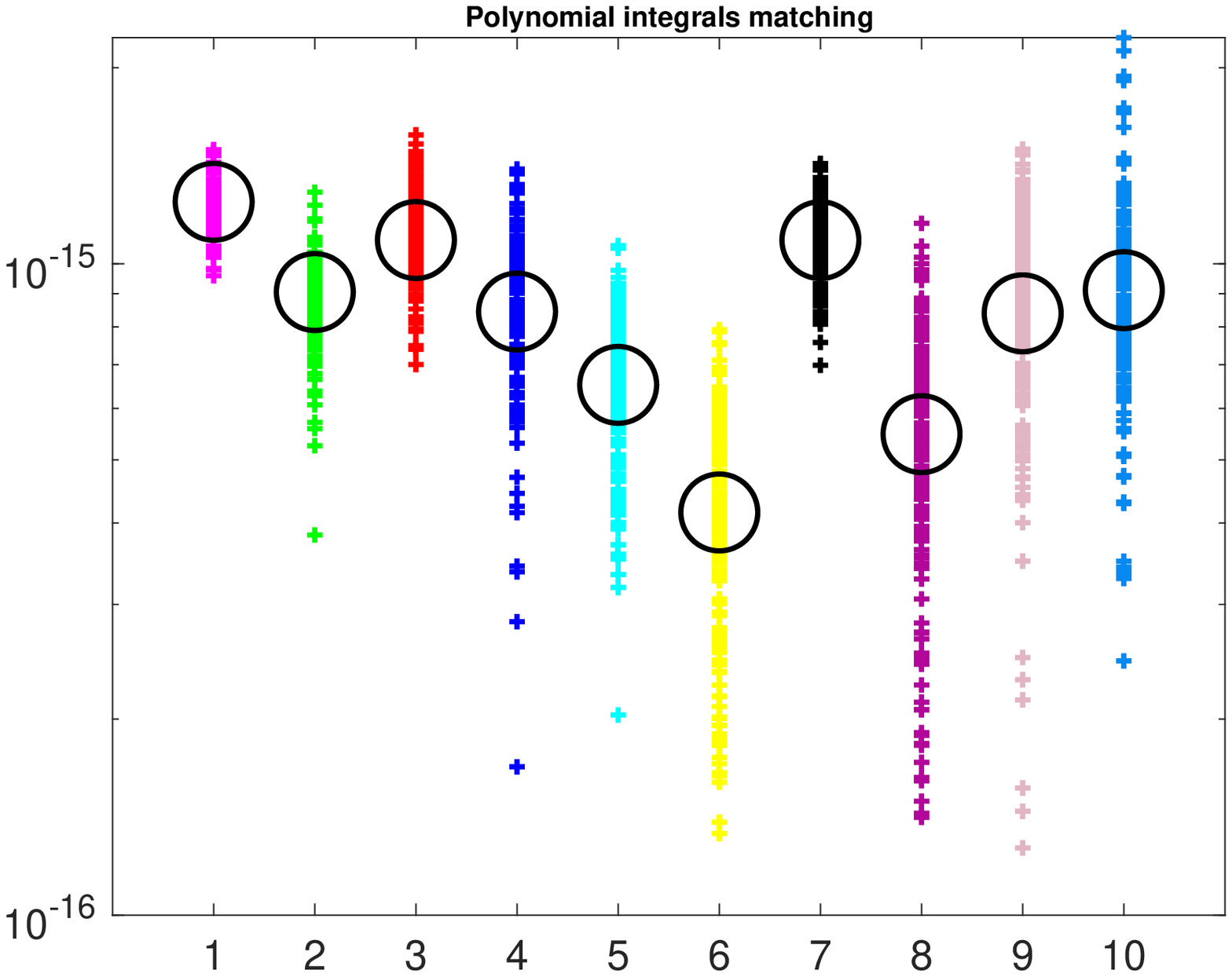}}
 \caption{Relative errors $E(g_k)$ of the tetrahedra-free rule over 100 polynomial integrands of the form $g_k=(a_kx+b_ky+c_kz+d_k)^n$ on the three polyhedra of Figure 1, where $a_k,b_k,c_k,d_k$ are uniform random coefficients in $[0,1]$ and $n=1,2,\dots,10$; the circles correspond to the average logarithmic error $\sum_{k=1}^{100}{\log(E(g_k))/100}$.}
 \label{fig_3CRI}
 \end{figure}

\begin{figure}[!htbp]
  \centering
   {\includegraphics[scale=0.32,clip]{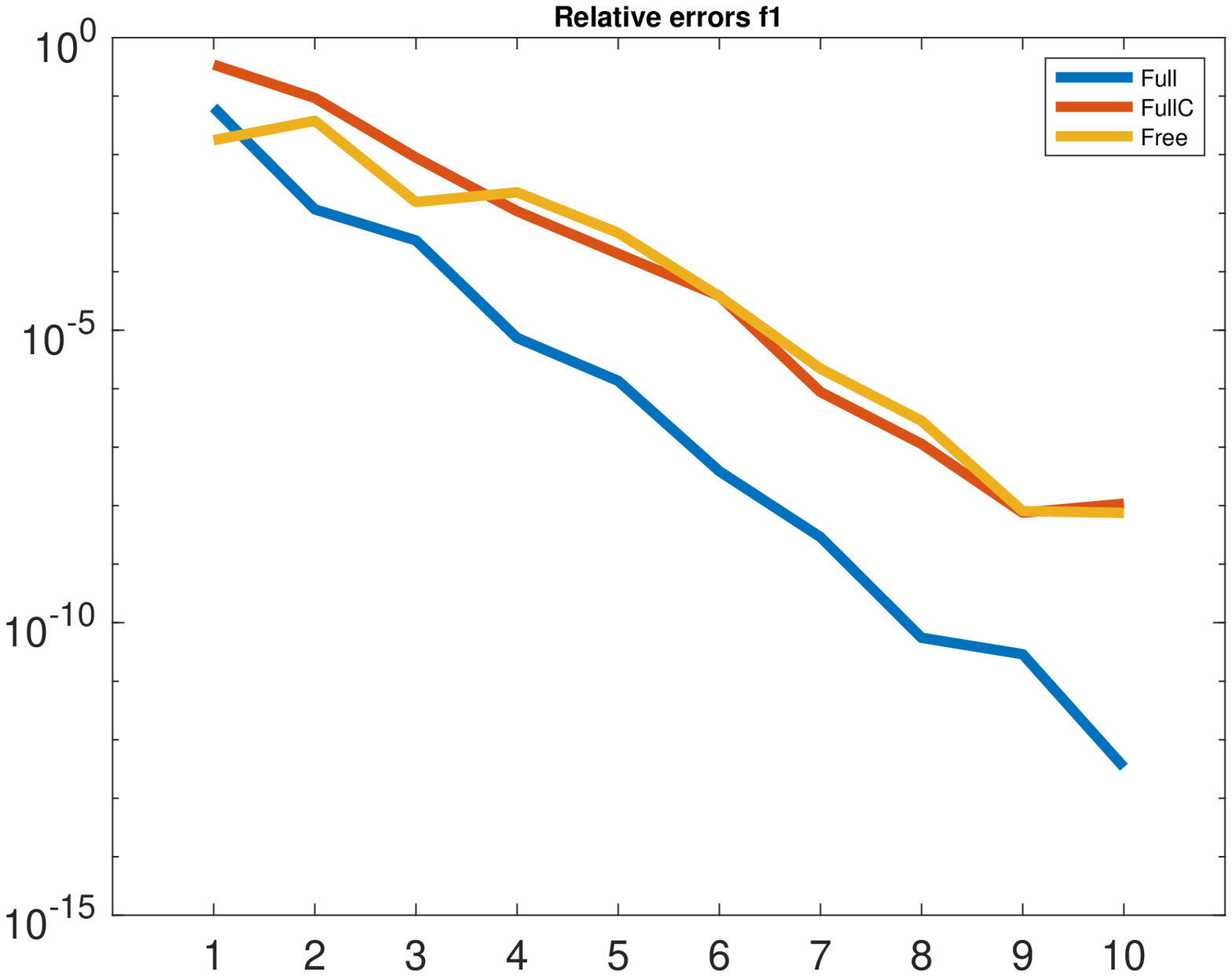}}
   {\includegraphics[scale=0.32,clip]{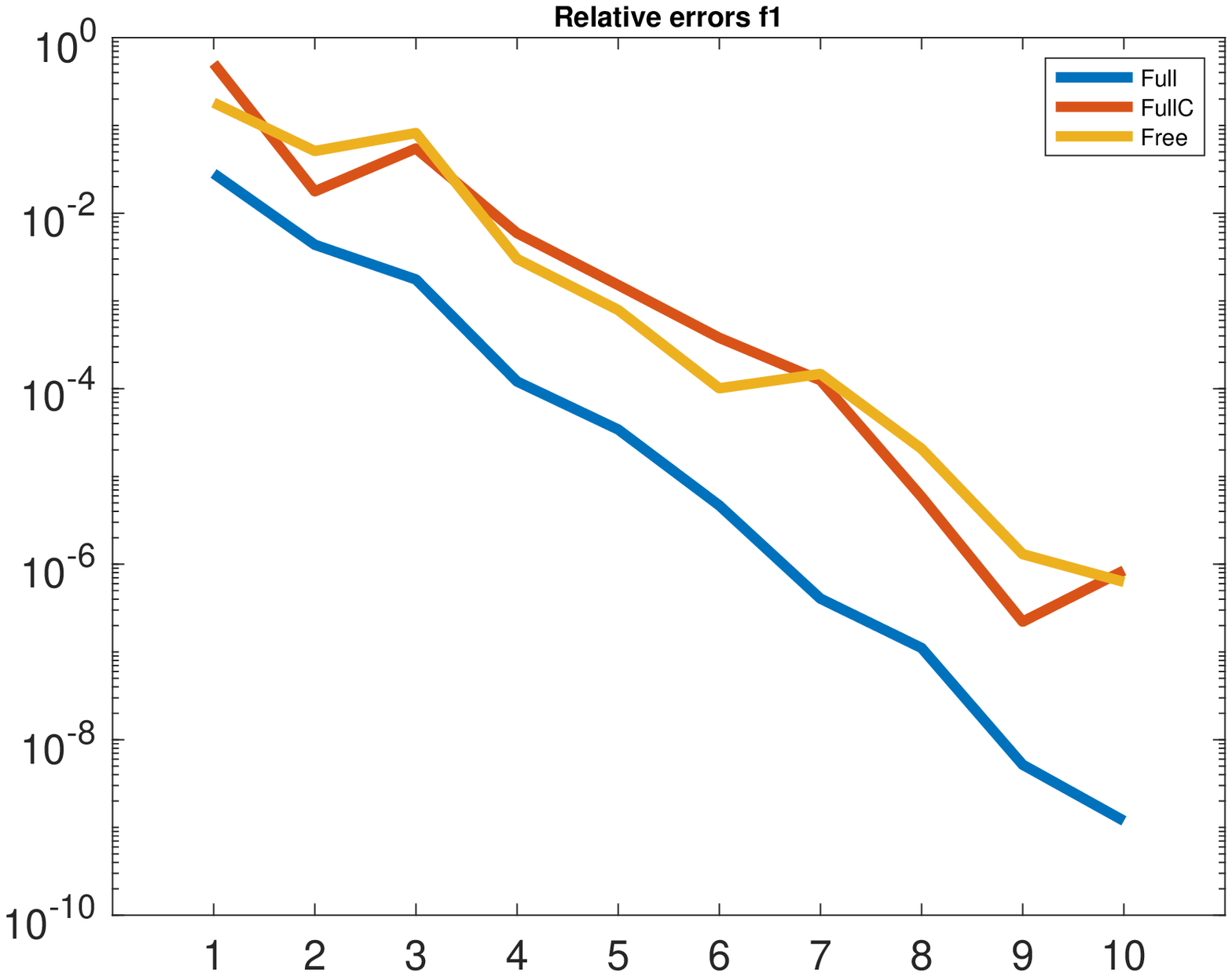}}
      {\includegraphics[scale=0.32,clip]{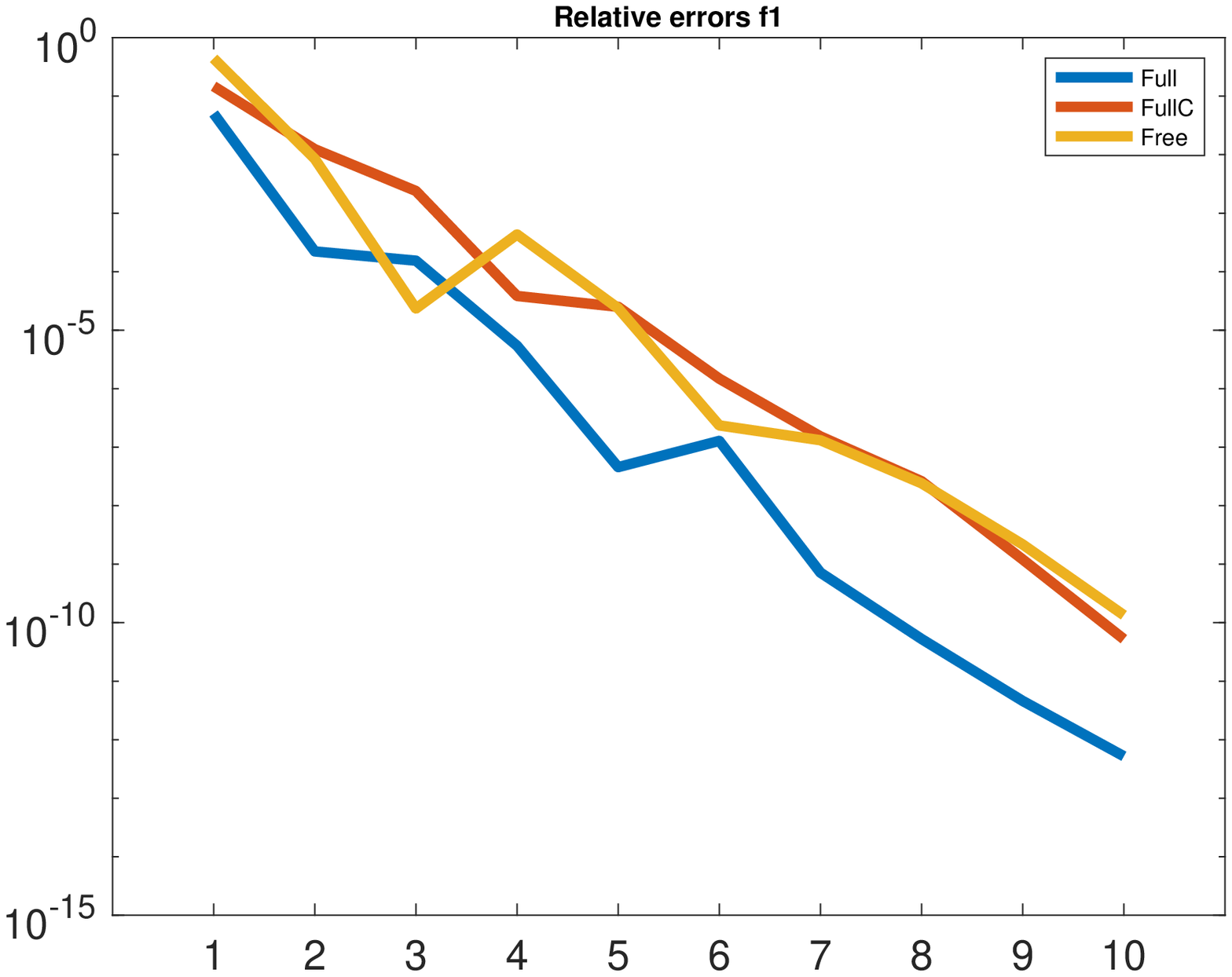}}
 \caption{Relative errors versus the exactness degree 
 in the integration of the Gaussian $f_1$ in (\ref{testf}) on the three polyhedra of Figure 1: uncompressed (blue) and compressed (brown) tetrahedra-based rules, tetrahedra-free rule(yellow).}
 \label{fig_3CRI}
 \end{figure}

 \begin{figure}[!htbp]
  \centering
   {\includegraphics[scale=0.32,clip]{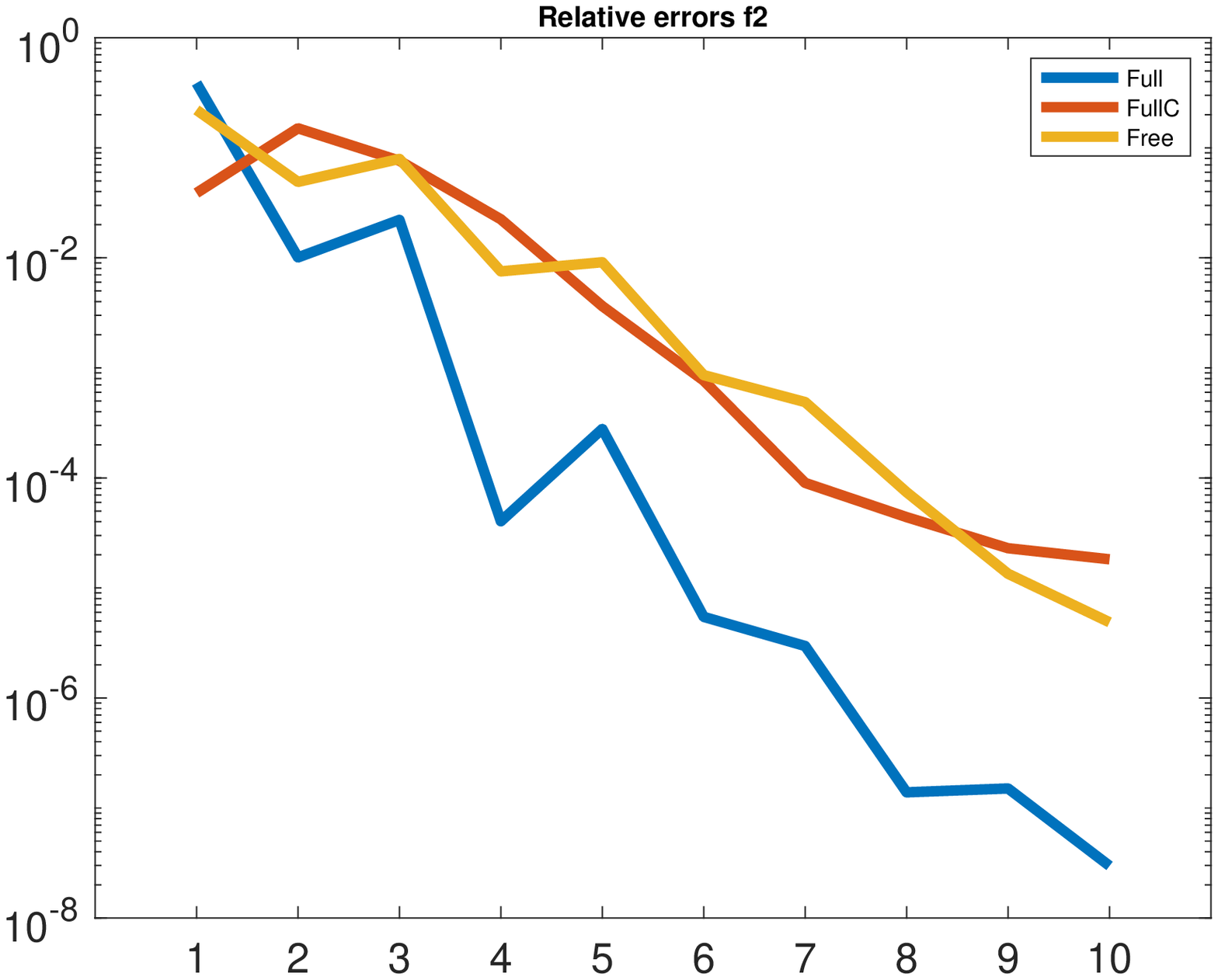}}
   {\includegraphics[scale=0.32,clip]{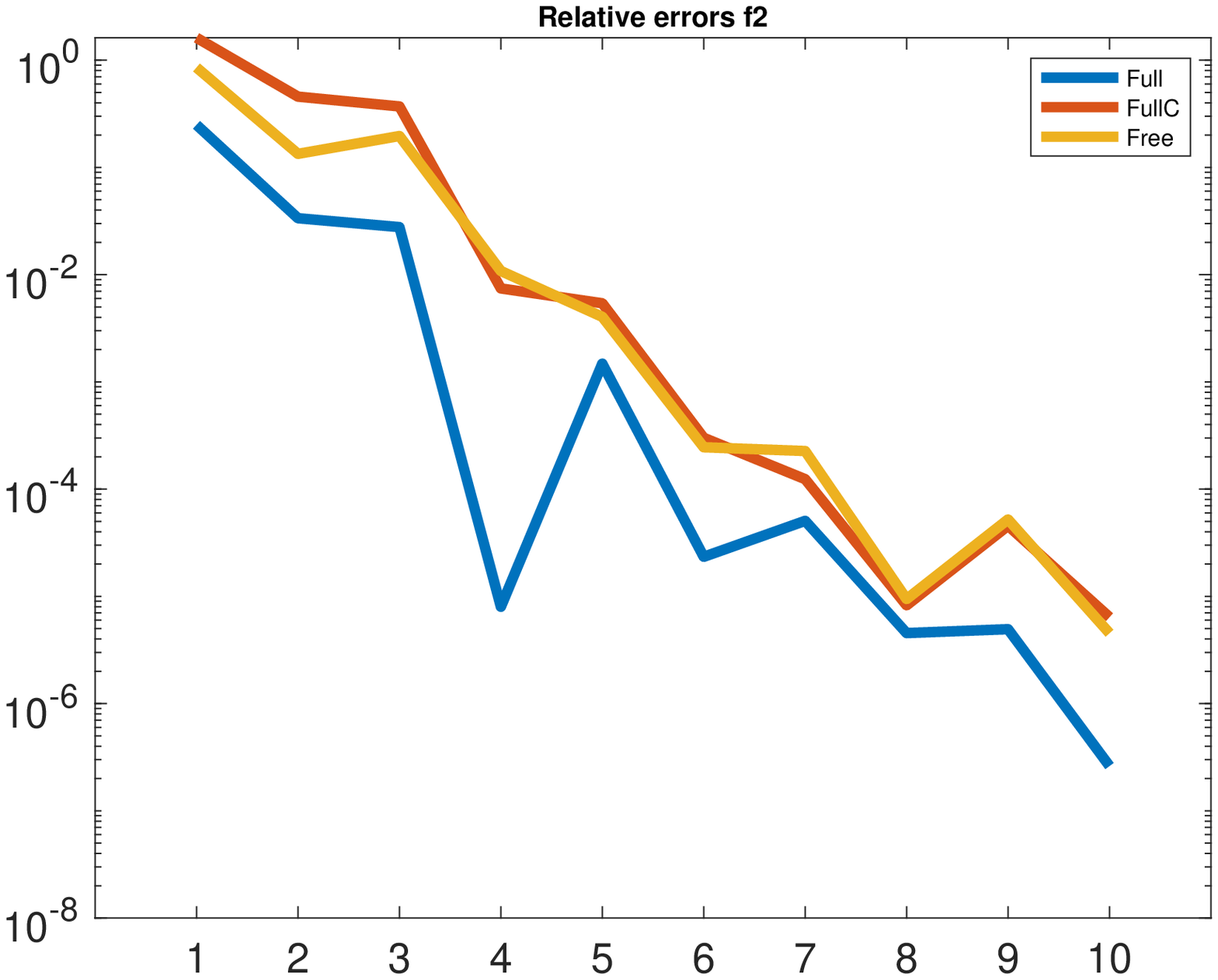}}
      {\includegraphics[scale=0.32,clip]{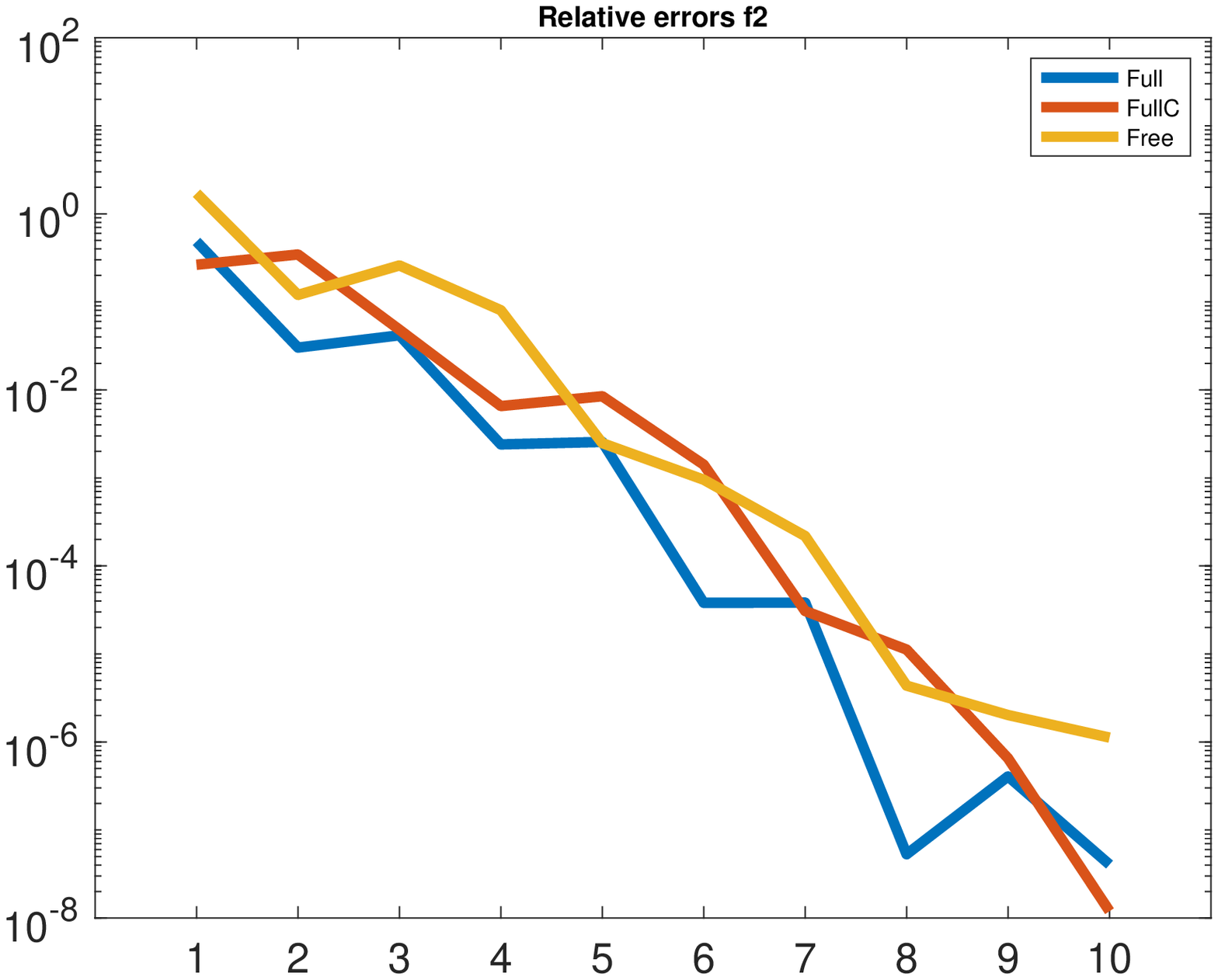}}
 \caption{Relative errors versus the exactness degree in the integration of the power function $f_2$ in (\ref{testf}) on the three polyhedra of Figure 1: uncompressed (blue) and compressed (brown) tetrahedra-based rules, tetrahedra-free rule(yellow).}
 \label{fig_3CRI}
 \end{figure}

 \begin{figure}[!htbp]
  \centering
   {\includegraphics[scale=0.32,clip]{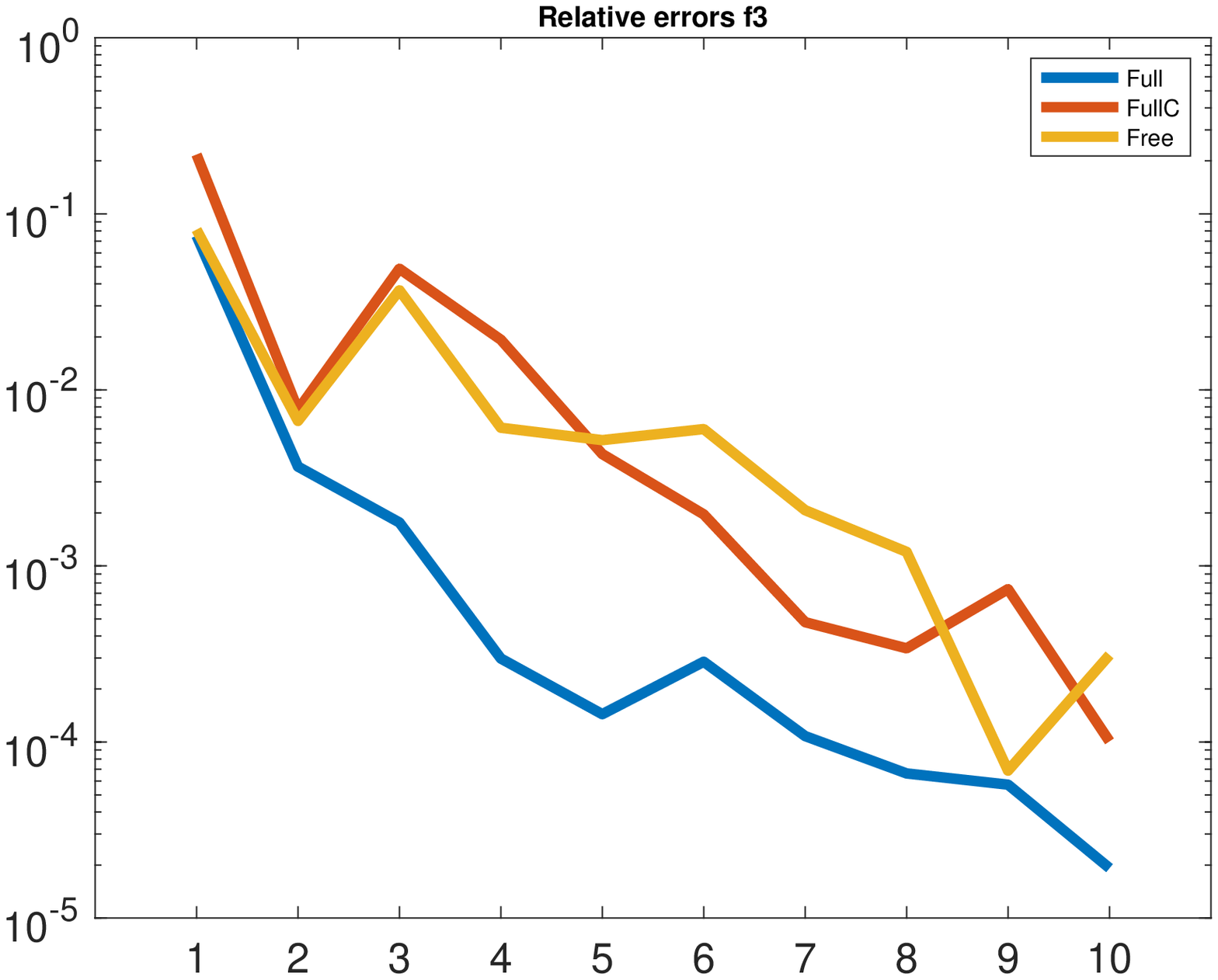}}
   {\includegraphics[scale=0.32,clip]{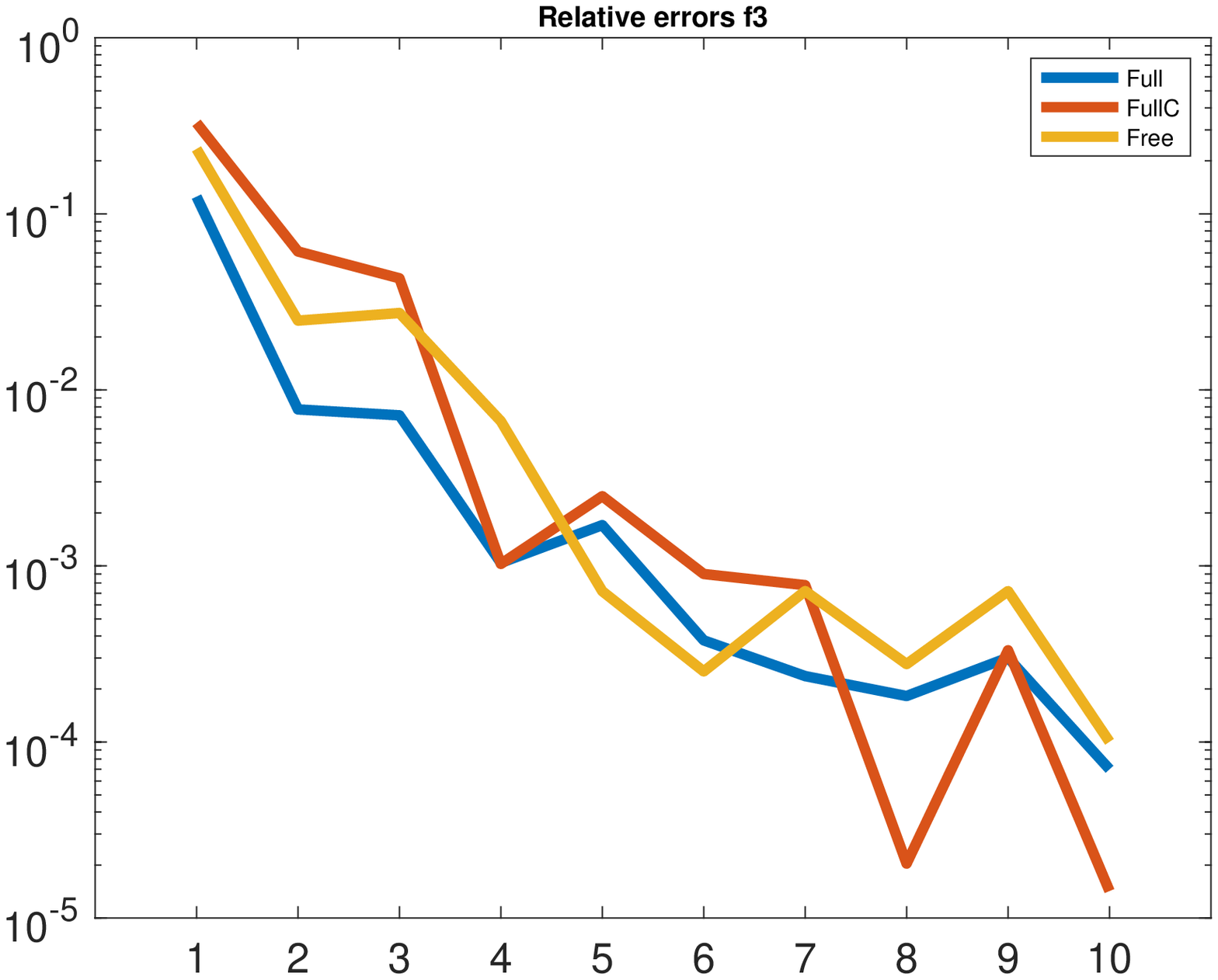}}
      {\includegraphics[scale=0.32,clip]{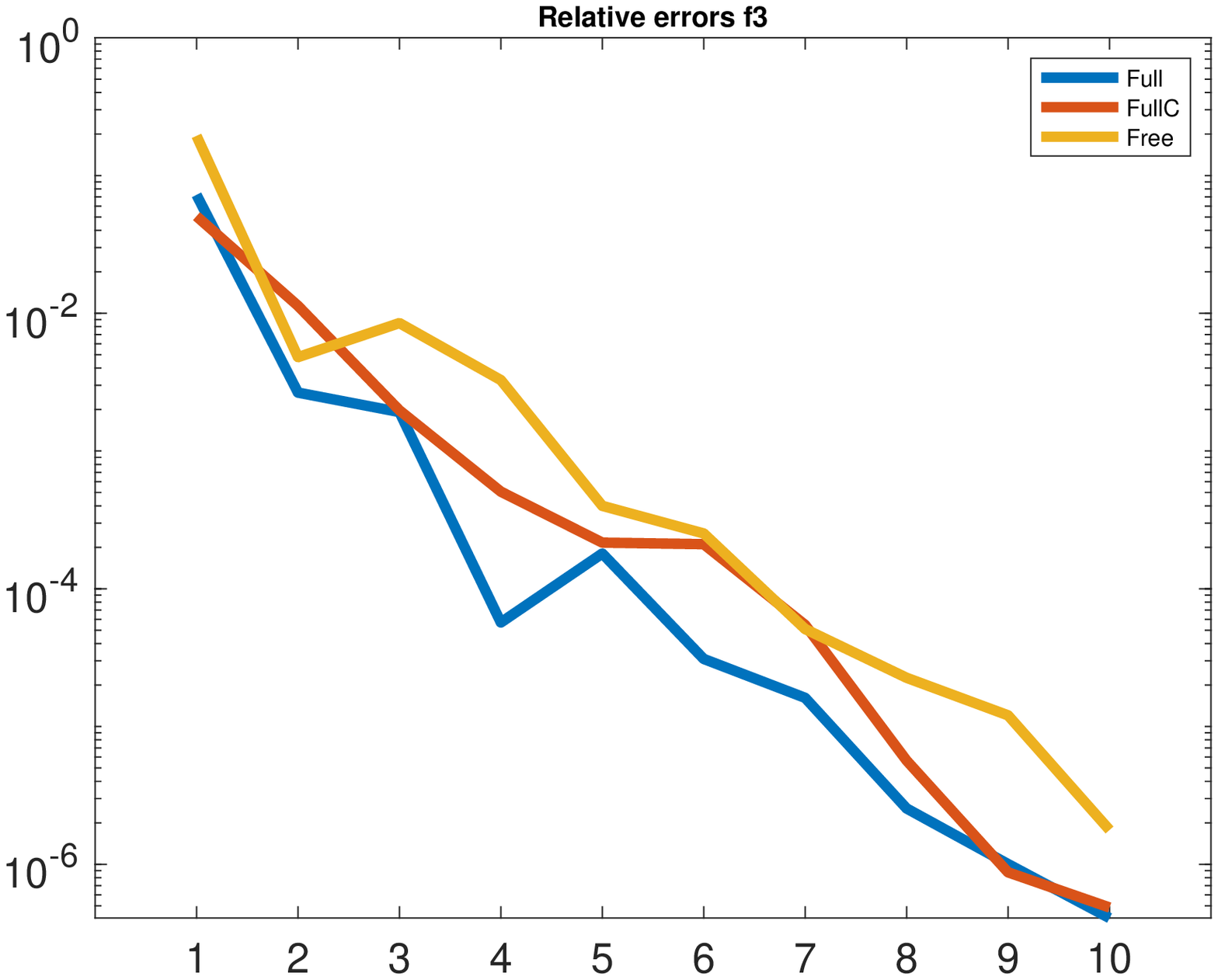}}
 \caption{Relative errors versus the exactness degree in the integration of the distance function $f_3$ in (\ref{testf}) on the three polyhedra of Figure 1: uncompressed (blue) and compressed (brown) tetrahedra-based rules, tetrahedra-free rule(yellow).}
 \label{fig_3CRI}
 \end{figure}

\vskip0.5cm
\begin{table}[ht]
\begin{center}\footnotesize
{\renewcommand{\arraystretch}{1.1}
\begin{tabular}{| c | c c c c c c c c c c | }
\hline
{\mbox{deg}} & $1$ & $2$ & $3$ & $4$ & $5$ & $6$ & $7$ & $8$ & $9$ & $10$    \\
\hline
$N$ & $4$ & $10$ & $20$ & $35$ & $56$ & $84$ & $120$ & $165$ & $220$ & $286$  \\
\hline
$ \Omega_1 $ & $31$ & $77$ & $160$ & $286$ & $447$ & $676$ & $969$ & $1327$ & $1765$ & $2286$  \\
$ \Omega_2 $ & $28$ & $78$ & $154$ & $273$ & $450$ & $669$ & $956$ & $1313$ & $1757$ & $2287$  \\
$ \Omega_3 $ & $29$ & $77$ & $156$ & $276$ & $446$ & $662$ & $950$ & $1312$ & $1753$ & $2285$  \\
 \hline
\end{tabular}
        }
\caption{\small{Cardinality of the Halton pointsets from which tetrahedra-free rules with $N$ nodes are obtained.}}
\label{tabT02}
\end{center}
\end{table}

 \section{Acknowledgements}

Work partially
supported by the
DOR funds and the biennial project BIRD 192932
of the University of Padova, and by the INdAM-GNCS 
2022 Project ``Methods and software for multivariate integral models''.
This research has been accomplished within the RITA ``Research ITalian network on Approximation", and within the UMI Group TAA ``Approximation Theory and Applications" (A. Sommariva).

\bibliographystyle{plain}

\end{document}